\documentclass[a4paper,10pt]{article}
\usepackage{amsmath}
\usepackage{amssymb}
\usepackage{amsthm}
\theoremstyle{plain}
\numberwithin{equation}{section}
\newtheorem{theo}{Theorem}[section]
\newtheorem{cor}{Corollary}[section]
\newtheorem{prop}{Proposition}[section]
\newtheorem{lem}{Lemma}[section]
\theoremstyle{definition}

\theoremstyle{remark}
\newtheorem{rem}{Remark}[section]
\begin{document}
\begin{Large}
{\bf Two hypotheses on the exponential class in 
the\\ class of $O$-subexponential 
infinitely divisible\\
 distributions}
\end{Large}
 \medskip
\begin{center} 
   Toshiro Watanabe\\
Center for Mathematical Sciences, The University of Aizu, Aizu-Wakamatsu, Fukushima 965-8580, Japan \\
E-mail: t-watanb@u-aizu.ac.jp   
   
 \end{center}
 {\bf Abstract}. Two hypotheses on the class $\mathcal{L}(\gamma)$ in the class $\mathcal{OS}\cap\mathcal{ID}$ are discussed.  Two weak hypotheses on the class $\mathcal{L}(\gamma)$  in the class $\mathcal{OS}\cap\mathcal{ID}$ are proved. A necessary and sufficient condition in order that, for  every $t>0$, the $t$-th convolution power of a distribution in the class $\mathcal{OS}\cap\mathcal{ID}$ belongs to  the class $\mathcal{L}(\gamma)$ is given. Sufficient conditions are given for the validity of two hypotheses on the class $\mathcal{L}(\gamma)$.\\
 
 
 
\section{Introduction and results}
\medskip
In what follows, we denote by $\mathbb R$ the real line and  by
 $\mathbb R_{+}$ the half line $[0,\infty)$. Denote by $\mathbb N$ the totality of positive integers and by $a\mathbb N$ the set $\{a,2a,3a,\dots \}$. The symbol $\delta_a(dx)$ stands for the delta measure at $a \in \mathbb R$. Let $\eta$ and $\rho$ be probability distributions on $\mathbb R$. We denote by $\eta*\rho$ the convolution of $\eta$ and $\rho$ and by $\rho^{n*}$
$n$-th convolution power of $\rho$ with the understanding that 
$\rho^{0*}(dx)=\delta_0(dx)$. Denote by $\bar\xi(x)$ the tail $\xi((x,\infty))$ of a measure $\xi$ on $\mathbb R$ for $x \in\mathbb R$. Let $\gamma \geq 0$. We define the $\gamma$-exponential moment $\widehat\xi(\gamma)$ as
$$\widehat\xi(\gamma):=\int_{-\infty}^{\infty}e^{\gamma x}\xi(dx).$$
If $\widehat\xi(\gamma) < \infty,$ we define the Fourier-Laplace transform $\widehat\xi(\gamma+iz)$ for $z\in \mathbb R$ as
$$\widehat\xi(\gamma+iz):=\int_{-\infty}^{\infty}e^{(\gamma + iz)x}\xi(dx).$$
An integral $\int_a^bg(x)\rho(dx)$ means $\int_{a+}^{b+}g(x)\rho(dx)$. For positive functions $f_1(x)$ and $g_1(x)$ on $[A,\infty)$ for some $A \in \mathbb R$, we define the relation $f_1(x)\sim g_1(x)$ by 
$\lim_{x \to \infty}f_1(x)/g_1(x)=1$ and the relation $f_1(x)\asymp g_1(x)$ by 
$$0<\liminf_{x \to \infty}f_1(x)/g_1(x)\leq \limsup_{x \to \infty}f_1(x)/g_1(x)< \infty.$$ 
Let $\gamma \geq 0$. A distribution $\rho$ on $\mathbb R$ belongs to the class $\mathcal{L}(\gamma)$ if $\overline{\rho}(x)>0$ for all $x >0$ and, for every $ a \in \mathbb R$,
$$\overline{\rho}(x+a)\sim e^{-\gamma a}\overline{\rho}(x).$$ 
 A distribution $\rho$ on $\mathbb R$ belongs to the class $\mathcal{S}(\gamma)$ if $\rho \in \mathcal{L}(\gamma)$, $ \widehat\rho(\gamma) < \infty$, and 
  $$\overline{\rho^{2*}}(x) \sim 2\widehat\rho(\gamma)\overline{\rho}(x).$$
  A distribution $\rho$ on $\mathbb R$ belongs to the class $\mathcal{OL}$ if $\overline{\rho}(x)>0$ for  $x >0$ and, for all $a \geq 0$,
$$\overline{\rho}(x-a) \asymp \overline{\rho}(x).$$
 A distribution $\rho$ on $\mathbb R$ belongs to the class $\mathcal{OS}$ if $\overline{\rho}(x)>0$ for all $x >0$ and  
 $$\overline{\rho^{2*}}(x) \asymp \overline{\rho}(x).$$
  Note that the class $\mathcal{OS}$ is included in the class  $\mathcal{OL}$.  A distribution $\rho$ on $\mathbb R$ belongs to the class $\mathcal{S}_{\sharp}$ if $\rho \in \mathcal{OS}$  and  
$$\limsup_{A \to \infty}\limsup_{x \to \infty}\frac{\overline{\rho}(x-A)\bar\rho(A)+\int_A^{x-A}\bar\rho(x-u)\rho(du)}{\bar\rho(x)}=0.$$
The class $\mathcal{S}_{\sharp}$ includes $\cup_{\gamma \geq 0}\mathcal{S}(\gamma)$ and it is   closed under convolution powers. 
A finite measure $\xi$ satisfies the {\it Wiener condition} if $\widehat\xi(iz)\ne 0$ for every $z \in \mathbb R$. Denote by $\mathcal{W}$ the totality of finite measures on $\mathbb R$ satisfying the Wiener condition.
 We denote by $\mathcal{ID}$ the class of all infinitely divisible distributions on $\mathbb R$. For $\mu \in \mathcal{ID}$, denote by $\nu$ its L\'evy measure. Under the assumption that $\bar\nu(c)>0$ for every $ c >0,$ define $\nu_1(dx):=1_{(1,\infty)}(x)\nu(dx)/\bar\nu(1).$ Let $\mu \in \mathcal{ID}$. We define a compound Poisson distribution $\mu_1$ with $c=\bar\nu(1)$ as
$$\mu_1(dx):=e^{-c}\sum_{k=0}^{\infty}\frac{c^k}{k!}\nu_1^{k*}(dx).$$ 
Denote by $\mu^{t*}$ the $t$-th convolution power of $\mu \in \mathcal{ID}$ for $t >0$. Note that $\mu^{t*}$ is the distribution of $X_t$ for a certain L\'evy process $\{X_t\}$ on $\mathbb R$. Let $\gamma \geq 0$. Define $T(\mu,\gamma)$ as
$$T(\mu,\gamma):=\{t>0: \mu^{t*} \in \mathcal{L}(\gamma)\}.$$
Since the class $\mathcal{L}(\gamma)$ is closed under convolutions by Theorem 3 of Embrechts and Goldie  \cite{eg1},  $T(\mu,\gamma)$ is empty or an additive semigroup in $(0,\infty)$.
We see from Lemma 2.2 below  that for $\mu \in \mathcal{OS}\cap\mathcal{ID}$, there are   positive integers $n$ such that  $\nu_1^{n*} \in \mathcal{OS}$. Let $n_0$ be the positive integer defined by (2.1) below. Note that we do not yet know an example of $\mu \in \mathcal{OS}\cap\mathcal{ID}$ such that $n_0 \ge 3.$

A class $\mathcal{C}$ of distributions is called {\it closed under convolution roots} if $\rho^{n*} \in \mathcal{C}$ for some $n \in \mathbb N$ implies $\rho \in \mathcal{C}$. We see from Shimura and Watanabe \cite{sw1} that the class $\mathcal{OS}$ is not  closed under convolution roots, but from Watanabe and Yamamuro \cite{wy1} that the class $\mathcal{OS}\cap\mathcal{ID}$ is closed under convolution roots. 
Embrechts et al.\ \cite{egv} in the one-sided case and Watanabe \cite{w1} in the two-sided case proved that the class $\mathcal{S}(0)$ is closed under convolution roots and Embrechts and Goldie \cite{eg1} conjectured that the class $\mathcal{L}(\gamma)$ with $\gamma \geq 0$ is   closed under convolution roots, but Shimura and Watanabe \cite{sw2} showed that the class $\mathcal{L}(\gamma)$ with $\gamma \geq 0$ is not  closed under convolution roots. Moreover, Watanabe and Yamamuro \cite{wy2} proved that the class $\mathcal{S}_{ac}$ of all absolutely  continuous distributions on $\mathbb R$ with subexponential densities  is not  closed under convolution roots. Embrechts and Goldie \cite{eg2} conjectured that the class $\mathcal{S}(\gamma)$ with $\gamma >0$ is   closed under convolution roots.  Watanabe \cite{w1} proved that $\mathcal{S}(\gamma)\cap\mathcal{ID}$ with $\gamma \geq 0$ is closed under convolution roots, but Watanabe \cite{w2}  showed that  the class $\mathcal{S}(\gamma)$ with $\gamma >0$ is not  closed under convolution roots. We add the following. Kl\"uppelberg \cite{k}  showed that the class $\mathcal{OS}$ is closed under convolutions. The class $\mathcal{S}(\gamma)$ is closed under convolution powers for $\gamma \geq 0$, but Leslie \cite{l}, for $\gamma = 0$, and  Kl\"uppelberg and Villasenor \cite{kv}, for $\gamma > 0$, proved that the class $\mathcal{S}(\gamma)$ is not closed under  convolutions. 

 We consider the following two  hypotheses on the class $\mathcal{L}(\gamma)$ in the class $\mathcal{OS}\cap\mathcal{ID}$ : 

\medskip
{\sc Hypothesis I}.  Let $\gamma \geq 0$. For every $\mu \in \mathcal{OS}\cap\mathcal{ID}$, if $\mu^{n*} \in \mathcal{L}(\gamma)$ for some $n \in \mathbb N$, then $\mu^{(n+1)*} \in \mathcal{L}(\gamma)$.

{\sc Hypothesis II}.   Let $\gamma \geq 0$. For every $\mu \in \mathcal{OS}\cap\mathcal{ID}$, if  $\mu^{n*} \in \mathcal{L}(\gamma)$ for some $n \in \mathbb N$, then $\mu \in \mathcal{L}(\gamma)$.

\medskip
We also consider the weak version of the above hypotheses :

\medskip
{\sc Hypothesis I'}.   Let $\gamma \geq 0$. For every $\mu \in \mathcal{OS}\cap\mathcal{ID}$, if $\mu^{n*}, \mu^{(n+1)*}  \in \mathcal{L}(\gamma)$ for some $n \in \mathbb N$, then $\mu^{(n+2)*} \in \mathcal{L}(\gamma)$.

{\sc Hypothesis II'}.   Let $\gamma \geq 0$. For every $\mu \in \mathcal{OS}\cap\mathcal{ID}$, if  $\mu^{n*}, \mu^{(n+1)*}  \in \mathcal{L}(\gamma)$ for some $n \in \mathbb N$, then $\mu \in \mathcal{L}(\gamma)$.

\medskip
Let $\gamma \geq 0$. Define
$$\mathcal{A}(\gamma):=\{\mu\in\mathcal{OS}\cap\mathcal{ID}:T(\mu,\gamma)=(0,\infty)\};$$
$$\mathcal{B}(\gamma):=\{\mu\in\mathcal{OS}\cap\mathcal{ID}:T(\mu,\gamma)=\emptyset\};$$
and
$$\mathcal{C}(\gamma):=\{\mu\in\mathcal{OS}\cap\mathcal{ID}:T(\mu,\gamma)=a_0\mathbb N \mbox{ with some }a_0 >0\}.$$

\begin{theo} Let $\gamma \geq 0$ and  $\mu \in \mathcal{OS}\cap\mathcal{ID}$. We have the following:

(i)  $\mathcal{OS}\cap\mathcal{ID}=\mathcal{A}(\gamma)\cup\mathcal{B}(\gamma)\cup\mathcal{C}(\gamma).$ Thus Hypotheses I' and II' are true.

(ii) The relation $\mu \in \mathcal{A}(\gamma)$ holds if and only if, for  all $a \geq 0$,
\begin{equation}
\lim_{x \to \infty}\frac{e^{-\gamma a}\overline{\nu_1}(x-a)-\overline{\nu_1}(x)}{\overline{\nu_1^{n_0*}}(x)} = 0.
\end{equation}
If $\mu \in \mathcal{A}(\gamma)$, then $\nu_1^{n*} \notin \mathcal{L}(\gamma)\cap\mathcal{OS}$ for $1 \leq n \leq n_0-1$ and $\nu_1^{n*} \in \mathcal{L}(\gamma)\cap\mathcal{OS}$ for $ n \geq n_0$.
\end{theo}
\begin{cor} Let $\gamma \geq 0$. Then the following are equivalent:

(1) Hypothesis I is true.

(2) Hypothesis II is true.

(3) $\mathcal{C}(\gamma)$ is empty.

(4) For every $\mu \in \mathcal{OS}\cap\mathcal{ID}$ it holds that, for every $2t \in T(\mu,\gamma)$ and for  every $a \geq 0$,
\begin{equation}
\limsup_{x \to \infty}\limsup_{\lambda \to \infty}\frac{|\int_x^{\lambda-x}(e^{-\gamma a}\overline{\mu_1^{t*}}(\lambda-a-u)-\overline{\mu_1^{t*}}(\lambda-u))\mu_1^{t*}(du)|}{\overline{\mu_1^{t*}}(\lambda)}=0.
\end{equation}
\end{cor}

\begin{rem} Let $\gamma = 0$. Then, $\mathcal{C}(0)$ is empty  and Hypotheses I and II are true. The relation $\mu \in \mathcal{A}(0)$ holds if and only if  
$$\lim_{x \to \infty}\frac{\nu_1((x,x+1])}{\overline{\nu_1^{n_0*}}(x)} = 0.$$
If $\mu \in \mathcal{A}(0)$, then $\nu_1^{n*} \notin \mathcal{L}(0)\cap\mathcal{OS}$ for $1 \leq n \leq n_0-1$ and $\nu_1^{n*} \in \mathcal{L}(0)\cap\mathcal{OS}$ for $ n \geq n_0$. Xu et al. showed in Theorem 2.2 of \cite{xfw} an example of $\mu \in \mathcal{A}(0)$ with $n_0=2$.
\end{rem}
For $\gamma > 0$, we cannot yet answer the question whether Hypotheses I and II are true.
    However, under some additional assumptions in terms of L\'evy measure, we establish that $\mathcal{C}(\gamma)$ is empty. 

\begin{prop} Let $\gamma > 0$ and $\mu \in \mathcal{OS}\cap\mathcal{ID}$.  Suppose that,  for every $a \geq 0$,
\begin{equation}
\liminf_{x\to\infty}e^{-\gamma a}\bar\nu_1(x-a)/\bar\nu_1(x) \geq 1.
\end{equation} 
Then, we have either $T(\mu,\gamma)= (0,\infty)$ or $\emptyset$. 

\end{prop} 
\begin{rem} Cui et al.\ \cite{cwx} proved a result analogous to the above proposition under a stronger assumption. Xu et al. showed in Theorem 1.1 of \cite{xwcy} an example of the case where $T(\mu,\gamma)\ne\emptyset$ in the above proposition.
\end{rem}
\begin{prop} Let $\gamma > 0$ and $\mu \in \mathcal{OS}\cap\mathcal{ID}$.  Suppose that $\nu_1^{2*} \in \mathcal{L}(\gamma)$ and the real part of $\widehat\nu_1(\gamma+iz)$ is not zero for every $z \in\mathbb R$. Then,  either $T(\mu,\gamma)= (0,\infty)$ or $\emptyset$. 

\end{prop}
\begin{prop} Let $\gamma > 0$ and $\mu \in \mathcal{OS}\cap\mathcal{ID}$.  Suppose that there exists $n_1 \in \mathbb N$ such that $\nu_1^{n_1*} \in \mathcal{S}_{\sharp}$. Then, either $T(\mu,\gamma)= (0,\infty)$ or $\emptyset$.  The equality   $T(\mu,\gamma)=(0,\infty)$ holds if and only if $\nu_1 \in \mathcal{S}(\gamma)$. 

\end{prop}
\begin{rem} Watanabe  made in Theorem 1.1 of \cite{w2} a distribution $\eta \in \mathcal{S}_{\sharp}$ such that $\eta^{n*} \in \mathcal{S}(\gamma)$ for every $n \geq 2$ but $\eta \notin \mathcal{S}(\gamma)$. Thus, taking  this $\eta$ as $\nu_1 $, then  Proposition 1.3 holds with $T(\mu,\gamma)= \emptyset$.
\end{rem}

\section{Preliminaries}
In this section, we give several basic results as preliminaries.  Pakes \cite{p1}
   proved the following.
\begin{lem}(Lemmas 2.1 and 2.5 of \cite{p1})
Let $\mu \in \mathcal{ID}$. Then we have
$\mu  \in \mathcal{L}(\gamma)$ if and only if $\mu_1 \in \mathcal{L}(\gamma)$.
\end{lem}
Watanabe and Yamamuro \cite{wy1} proved the following.
\begin{lem}(Proposition 3.1 of \cite{wy1}) Suppose that $\mu \in \mathcal{ID}$. Then, we have
$\mu \in \mathcal{OS}$ if and only if there is $n \in \mathbb N$ such that $\nu_1^{n*} \in \mathcal{OS}$ and $\overline{\mu_1^{t*}}(x) \asymp \overline{\nu_1^{n*}}(x)$ for any $t >0$.
\end{lem}
For $\mu \in \mathcal{OS}\cap\mathcal{ID}$, define $n_0 \in \mathbb N$
as 
\begin{equation}
n_0:=\min\{n \in \mathbb N: \nu_1^{n*} \in \mathcal{OS}\}.
\end{equation}
\begin{lem} Let $\mu \in \mathcal{OS}\cap\mathcal{ID}$. 

(i) There exists $C(a) >0$ such that, for all $ a \geq 0$ and all $ x >0$,
$$\overline{\nu_1^{n_0*} }(x-a)\leq C(a) \overline{\nu_1^{n_0*} }(x).$$

(ii) There exists $K >1$ such that, for all $n \in \mathbb N$ and all $ x >0$,
$$\overline{\nu_1^{n*} }(x)\leq K^n \overline{\nu_1^{n_0*} }(x).$$
\end{lem}
Proof. Assertion (i) is clear since $\nu_1^{n_0*} \in \mathcal{OS}\subset\mathcal{OL}.$  We see from Proposition 2.4 of Shimura and Watanabe \cite{sw1} that there exists $K_1 >1$ such that, for all $k \in \mathbb N$ and all $ x >0$,
$$\overline{\nu_1^{(kn_0)*} }(x)\leq K_1^k \overline{\nu_1^{n_0*} }(x).$$
Note that, for $m \leq n$, 
$$\overline{\nu_1^{m*} }(x)\leq  \overline{\nu_1^{n*} }(x).$$
Hence, we have, for $0 \leq j \leq n_0-1$ and for all $k \in \mathbb N$, with $K=K_1^{2/n_0}>1$
$$\overline{\nu_1^{(kn_0+j)*} }(x)\leq K^{(kn_0+j)} \overline{\nu_1^{n_0*} }(x).$$
This inequality holds for $k=0$ too. Thus assertion (ii) is true. \hfill $\Box$

Under the assumption that $\zeta \in\mathcal{OS}\subset \mathcal{OL}$, we define the following. Let 
$$d^*:=\limsup_{x \to \infty}
\frac{\overline{\zeta^{2*}}(x)}{\overline{\zeta}(x)}<\infty .$$
Let $\Lambda$ be the totality of increasing sequences $\{\lambda_n\}_{n=1}^{\infty}$ with $\lim_{n \to\infty} \lambda_n = \infty$ such that, for every $x \in \mathbb R$, the following limit exists and is finite:
\begin{equation}
m(x;\{\lambda_n\}):=\lim_{n \to\infty}\frac{\bar\zeta(\lambda_n-x)} {\bar\zeta(\lambda_n)}.
\end{equation}
Define, for each sequence   
 $\{x_n\}_{n=1}^{\infty}$ with $\lim_{n \to\infty} x_n = \infty$, $T_n(y)$ as
$$ T_n(y):=\frac{\bar\zeta(x_n-y)} {\bar\zeta(x_n)}.$$
Since $\{T_n(y)\}_{n=1}^{\infty}$ is a sequence of increasing functions, uniformly bounded on every finite interval, by Helly's  selection principle,  there exists an increasing subsequence $\{\lambda_n\}$ of $\{x_n\}$ with $\lim_{n \to\infty} \lambda_n = \infty$ such that everywhere on $\mathbb R$ (2.2) holds.
 The limit function $m(x;\{\lambda_n\})$ is increasing and is  finite. That is, $\{\lambda_n\} \in \Lambda$. It follows  that, under the assumption that $\zeta \in\mathcal{OS}$, there exists an increasing subsequence $\{\lambda_n\} \in \Lambda $ of $\{x_n\}$ for each sequence   
 $\{x_n\}_{n=1}^{\infty}$ with $\lim_{n \to\infty} x_n = \infty$.
\begin{lem}
Suppose that $\zeta \in\mathcal{OS}$. Then, we have the following.

(i) If $\{\lambda_n\} \in \Lambda $, then $\{\lambda_n-a\} \in \Lambda $ for every $a \in \mathbb R$.

(ii) For  $\{\lambda_n\} \in \Lambda $,
$$\int_{-\infty}^{\infty}m(x;\{\lambda_n\})\zeta(dx) < \infty$$
and
$$\lim_{a \to \infty}m(a;\{\lambda_n\})\bar\zeta(a)=0.$$
In particular, if $\zeta \in\mathcal{OS}\cap\mathcal{L}(\gamma)$, then 
$m(x;\{\lambda_n\})=e^{\gamma x}$ and $\widehat\zeta(\gamma) < \infty.$
\end{lem}
Proof. We prove (i). Suppose that $\{\lambda_n\} \in \Lambda. $ We have, for $x,a \in \mathbb R$,
$$\lim_{n \to\infty}\frac{\bar\zeta(\lambda_n-a-x)} {\bar\zeta(\lambda_n-a)}=\frac{m(x+a;\{\lambda_n\})}{m(a;\{\lambda_n\})}.$$
Thus $\{\lambda_n-a\} \in \Lambda. $
Next, we prove (ii). Let $\rho$ be a distribution on $\mathbb R$. Note that, for $x >2A,$
\begin{equation}
\overline{\rho^{2*}}(x)=2\int_{-\infty}^{A+}\bar\rho(x-u)\rho(du)+\overline{\rho}(x-A)\bar\rho(A)+\int_A^{x-A}\bar\rho(x-u)\rho(du).
\end{equation}
We see from (2.3) that, for  $\{\lambda_n\} \in \Lambda $ and $s >0$,
\begin{equation}
\begin{split}
d^* &\geq\limsup_{n \to \infty}
\frac{\overline{\zeta^{2*}}(\lambda_n)}{\overline{\zeta}(\lambda_n)}\\
\nonumber
&\geq 2\limsup_{n \to \infty}\int_{-\infty}^{s+}\frac{\bar\zeta(\lambda_n-x)} {\bar\zeta(\lambda_n)}\zeta(dx)\\
&\geq 2\int_{-\infty}^{s+}m(x;\{\lambda_n\})\zeta(dx).
\end{split}
\end{equation}
As $ s \to \infty$, we have 
$$\int_{-\infty}^{\infty}m(x;\{\lambda_n\})\zeta(dx) < \infty.$$
Since $m(x;\{\lambda_n\})$ is increasing in $x$, we have
\begin{equation}
\begin{split}
&\lim_{a \to \infty}m(a;\{\lambda_n\})\bar\zeta(a)\\ \nonumber
& \leq\lim_{a \to \infty}\int_{a+}^{\infty}m(x;\{\lambda_n\})\zeta(dx)=0.
\end{split}
\end{equation}
Hence, if $\zeta \in\mathcal{OS}\cap\mathcal{L}(\gamma)$, then 
$m(x;\{\lambda_n\})=e^{\gamma x}$ and $\widehat\zeta(\gamma) < \infty.$
Thus we have proved the lemma.     \hfill $\Box$

Pakes \cite{p1,p2} asserted and Watanabe \cite{w1} finally proved the following.
\begin{lem}(Theorem 1.1 of \cite{w1}) Let $\gamma \geq 0$.  Then $\mu \in \mathcal{ID}\cap\mathcal{S}(\gamma)$ if and only if $\nu_1 \in \mathcal{S}(\gamma)$.
\end{lem}

\begin{lem} Let $\gamma \geq 0$.  Suppose that $\rho \in \mathcal{S}_{\sharp}$. 

(i) If $\bar\eta(x) \asymp \bar\rho(x)$, then $\eta \in \mathcal{S}_{\sharp}$.

(ii)  $\rho \in \mathcal{S}(\gamma)$ if and only if $\rho \in \mathcal{L}(\gamma)$.
\end{lem}
Proof.  Suppose that $\rho \in \mathcal{S}_{\sharp}$. We prove (i). If $\bar\eta(x) \asymp \bar\rho(x)$, then there is $C >0$ such that  $\bar\eta(x)\leq C \bar\rho(x)$ for $x \in\mathbb R$. By using integration by parts in the second inequality, we obtain that
\begin{equation}
\begin{split}
&\bar\eta(x-A)\bar\eta(A)+\int_A^{x-A}\bar\eta(x-u)\eta(du)\\ \nonumber
&\leq C^2\bar\rho(x-A)\bar\rho(A)+C \int_A^{x-A}\bar\rho(x-u)\eta(du)\\
&\leq 2C^2\bar\rho(x-A)\bar\rho(A)+C^2\int_A^{x-A}\bar\rho(x-u)\rho(du).
\end{split}
\end{equation}
Thus, we see that
$$\limsup_{A \to \infty}\limsup_{x \to \infty}\frac{(\overline{\eta}(x-A)\bar\eta(A)+\int_A^{x-A}\bar\eta(x-u)\eta(du))}{\bar\eta(x)}=0.$$
That is, $\eta \in  \mathcal{S}_{\sharp}$.
Next we prove (ii).
If $\rho \in \mathcal{S}(\gamma)$, then clearly $\rho \in \mathcal{L}(\gamma)$. Note that, for $x >2A,$ (2.3) holds.
If $\rho \in \mathcal{S}_{\sharp}\cap \mathcal{L}(\gamma)$, then we have
\begin{equation}
\begin{split}
&\lim_{x \to \infty}\frac{\overline{\rho^{2*}}(x)}{\bar\rho(x)}\\
\nonumber
=&\lim_{A \to\infty}2\int_{-\infty}^{A+}\lim_{x \to \infty}\frac{\bar\rho(x-u)}{\bar\rho(x)}\rho(du)\\
=&2\widehat\rho(\gamma) < \infty.
\end{split}
\end{equation}
Thus we see that $\rho \in \mathcal{S}(\gamma)$.\hfill $\Box$

Watanabe \cite{w2} extended Wiener's approximation theorem in \cite{w} as follows.
\begin{lem} (Lemma 2.6 of Watanabe \cite{w2})
Let $\xi$ be a finite measure on $\mathbb R$. The following are equivalent:

(1)  $\xi \in \mathcal{W}.$

(2)  If, for a bounded measurable function $g(x)$ on $\mathbb R$,
$$\int_{-\infty}^{\infty}g(x-t)\xi(dt)=0 \mbox{ for a.e. } x\in\mathbb R,$$
then $g(x)=0$  for a.e. $x\in\mathbb R.$
\end{lem}
\section{Convolution lemmas}
In this section, we give important lemmas on convolutions.
\begin{lem} Let $\gamma \geq 0$. Suppose that $\zeta \in \mathcal{OS} $. For $j=1,2$, let $\rho_j$  be distributions on $\mathbb R_+$  satisfying
\begin{equation}
 \bar\rho_j(x) \leq C_j\bar\zeta(x)\mbox{ with some
 } C_j>0 \mbox{ for all  } x >0.
\end{equation}
 Let $\{\lambda_n\} \in \Lambda. $ 
 
 (i) Let $\lambda_n > a+x$ and $x >0$.  We have, for every $a \geq 0$,
\begin{equation}
e^{-\gamma a}\overline{\rho_1*\rho_2}(\lambda_n-a)-\overline{\rho_1*\rho_2}(\lambda_n)=:\sum_{j=1}^{4}I_j,
\end{equation}
where
$$I_1:=-\int_{\lambda_n-a-x}^{\lambda_n-x}\overline{\rho_1}(\lambda_n-y)\rho_2(dy),$$
$$I_2:=\overline{\rho_1}(x)(e^{-\gamma a}\overline{\rho_2}(\lambda_n-a-x)-\overline{\rho_2}(\lambda_n-x)),$$
$$I_3:=\int_{0-}^{(\lambda_n-a-x)+}(e^{-\gamma a}\overline{\rho_1}(\lambda_n-a-y)-\overline{\rho_1}(\lambda_n-y))\rho_2(dy),$$
and
$$I_4:=\int_{0-}^{x+}(e^{-\gamma a}\overline{\rho_2}(\lambda_n-a-y)-\overline{\rho_2}(\lambda_n-y))\rho_1(dy).$$

(ii) We have for $j=1,2$
\begin{equation}
\limsup_{x\to \infty}\limsup_{n \to\infty}\frac{|I_j|}{\bar\zeta(\lambda_n)}=0.
\end{equation}

\end{lem}
Proof. By using integration by parts, we have
\begin{equation}
\begin{split} 
&\overline{\rho_1*\rho_2}(\lambda_n-a)\\ \nonumber
&=\int_{0-}^{(\lambda_n-a-x)+}\overline{\rho_1}(\lambda_n-a-y)\rho_2(dy)\\
&+\int_{\lambda_n-a-x}^{\lambda_n-a}\overline{\rho_1}(\lambda_n-a-y)\rho_2(dy)+\overline{\rho_2}(\lambda_n-a)\\
&=\int_{0-}^{(\lambda_n-a-x)+}\overline{\rho_1}(\lambda_n-a-y)\rho_2(dy)+\int_{0-}^{x+}\overline{\rho_2}(\lambda_n-a-y)\rho_1(dy)\\
&+\overline{\rho_1}(x)\overline{\rho_2}(\lambda_n-a-x),
\end{split}
\end{equation} 
and
\begin{equation}
\begin{split} 
&\overline{\rho_1*\rho_2}(\lambda_n)\\ \nonumber
&=\int_{0-}^{(\lambda_n-a-x)+}\overline{\rho_1}(\lambda_n-y)\rho_2(dy)
+\int_{\lambda_n-a-x}^{\lambda_n-x}\overline{\rho_1}(\lambda_n-y)\rho_2(dy)\\
&+\int_{\lambda_n-x}^{\lambda_n}\overline{\rho_1}(\lambda_n-y)\rho_2(dy)
+\overline{\rho_2}(\lambda_n)\\
&=\int_{0-}^{(\lambda_n-a-x)+}\overline{\rho_1}(\lambda_n-y)\rho_2(dy)+\int_{\lambda_n-a-x}^{\lambda_n-x}\overline{\rho_1}(\lambda_n-y)\rho_2(dy)\\
&+\int_{0-}^{x+}\overline{\rho_2}(\lambda_n-y)\rho_1(dy)+\overline{\rho_1}(x)\overline{\rho_2}(\lambda_n-x).
\end{split}
\end{equation}
Thus assertion (i) is valid. We have by Lemma 2.4 for $j=1,2$
\begin{equation}
\begin{split}
&\limsup_{x\to \infty}\limsup_{n \to\infty}\frac{|I_j|}{\bar\zeta(\lambda_n)}\\ \nonumber
&\leq \limsup_{x \to\infty}\overline{\rho_1}(x)\limsup_{n \to \infty}\frac{\overline{\rho_2}(\lambda_n-a-x)}{\bar\zeta(\lambda_n)}\\
&\leq C_1C_2\limsup_{x \to\infty}\bar\zeta(x)m(x;\{\lambda_n-a\})m(a;\{\lambda_n\})=0.
\end{split}
\end{equation}

 \begin{lem} Let $\gamma \geq 0$.  Suppose that $\zeta \in \mathcal{OS} $. For $j=1,2$, let $\rho_j$  be distributions on $\mathbb R_+$  satisfying (3.1).
Suppose further that for $j=1,2$ and  every $a \geq 0$,
$$\lim_{x \to \infty}\frac{e^{-\gamma a}\bar\rho_j(x-a)-\bar\rho_j(x)}{\bar\zeta(x)} = 0.$$
Then, for every $a \geq 0$,
\begin{equation}
\lim_{x \to \infty}\frac{e^{-\gamma a}\overline{\rho_1*\rho_2}(x-a)-
\overline{\rho_1*\rho_2}(x)}{\bar\zeta(x)} = 0.
\end{equation}
\end{lem}
Proof. Let $\{\lambda_n\} \in \Lambda. $ By the assumption for $j=1$,
 there is $\epsilon(x)>0 $ such that $\epsilon(x) \to 0$ as $x \to \infty$ and
$$|e^{-\gamma a}\overline{\rho_1}(\lambda_n-a-y)-\overline{\rho_1}(\lambda_n-y)|\leq\epsilon(x)\zeta(\lambda_n-y)$$
for $0 \leq y \leq\lambda_n-a-x.$
Thus we have
\begin{equation}
\begin{split}
&\limsup_{x \to\infty}\limsup_{n \to \infty}\frac{|I_3|}{\bar\zeta(\lambda_n)}\\
&\leq \limsup_{x \to\infty}\epsilon(x) \limsup_{n \to \infty}\frac{\int_{0-}^{(\lambda_n-a-x)+}\bar\zeta(\lambda_n-y)\rho_2(dy)}{\bar\zeta(\lambda_n)}\\
&\leq\limsup_{x \to\infty}\epsilon(x)\limsup_{n \to \infty}\frac{ \bar\zeta(\lambda_n)+\int_{a+x}^{\lambda_n}\overline{\rho_2}(\lambda_n-y)\zeta(dy)}{\bar\zeta(\lambda_n)}\\
&\leq\limsup_{x \to\infty}\epsilon(x) \limsup_{n \to \infty}\frac{\bar\zeta(\lambda_n)+C_2\overline{\zeta^{2*}}(\lambda_n)}{\bar\zeta(\lambda_n)}=0.
\end{split}
\end{equation}
As in the above argument, we have
\begin{equation}
\limsup_{x \to\infty}\limsup_{n \to \infty}\frac{|I_4|}{\bar\zeta(\lambda_n)}=0.\nonumber
\end{equation}
Thus, by  (3.2) and (3.3) of Lemma 3.1, we have proved (3.4). \hfill $\Box$
\begin{lem} Let $\gamma \geq 0$.  Suppose that $\zeta \in \mathcal{OS} $. For $j=1,2$, let $\rho_j$  be distributions on $\mathbb R_+$  satisfying (3.1). Suppose further that, for $j=1,2,$ and for every $a \geq 0$,
$$\liminf_{x \to \infty}\frac{e^{-\gamma a}\bar\rho_j(x-a)-\bar\rho_j(x)}{\bar\zeta(x)} \geq 0.$$
Then, we have, for every $a \geq 0$,
\begin{equation}
\liminf_{x \to \infty}\frac{e^{-\gamma a}\overline{\rho_1*\rho_2}(x-a)-\overline{\rho_1*\rho_2}(x)}{\bar\zeta(x)} \geq 0.
\end{equation}
\end{lem}
Proof. Let  $\{\lambda_n\} \in \Lambda. $  Let $\epsilon >0$ and  $a \geq 0$ be arbitrary and let $ n \in \mathbb N$ and $x \in (0,\lambda_n-a)$ be sufficiently large such that 
$$e^{-\gamma a}\overline{\rho_1}(\lambda_n-a-y)-\overline{\rho_1}(\lambda_n-y) \geq -\epsilon \bar\zeta(\lambda_n-y)$$
for $0 \leq y \leq \lambda_n-a-x$ and
$$e^{-\gamma a}\overline{\rho_2}(\lambda_n-a-y)-\overline{\rho_2}(\lambda_n-y)\geq -\epsilon\bar\zeta(\lambda_n-y)$$
for $0 \leq y \leq x.$
By (3.2) and (3.3) of Lemma 3.1, we have only to prove that
$$\sum_{j=3}^{4}\liminf_{x \to\infty}\liminf_{n \to \infty}\frac{I_j}{\bar\zeta(\lambda_n)}
\geq 0.$$
We have
\begin{equation}
\begin{split}
I_3&\geq -\epsilon \int_{0-}^{(\lambda_n-a-x)+}\bar\zeta(\lambda_n-y)\rho_2(dy)\\ \nonumber
&\geq-\epsilon \left(\bar\zeta(\lambda_n)+\int_{a+x}^{\lambda_n}\overline{\rho_2}(\lambda_n-y)\zeta(dy)\right)\\
&\geq-\epsilon \left(\bar\zeta(\lambda_n)+C_2\overline{\zeta^{2*}}(\lambda_n)\right),
\end{split}
\end{equation}
and
\begin{equation}
\begin{split}
I_4&\geq -\epsilon \int_{0-}^{x+}\bar\zeta(\lambda_n-y)\rho_1(dy)\\
\nonumber
&\geq-\epsilon \left(\bar\zeta(\lambda_n)+\int_{\lambda_n-x}^{\lambda_n}\overline{\rho_1}(\lambda_n-y)\zeta(dy)\right)\\
&\geq-\epsilon \left(\bar\zeta(\lambda_n)+C_1\overline{\zeta^{2*}}(\lambda_n)\right).
\end{split}
\end{equation}
Thus we see that
$$\liminf_{n \to \infty}\frac{I_3}{\bar\zeta(\lambda_n)}\geq-\epsilon (1+C_2d^*),$$
and
$$\liminf_{n \to \infty}\frac{I_4}{\bar\zeta(\lambda_n)}\geq-\epsilon(1+C_1d^*).$$ 
Since $\epsilon >0$ is arbitrary, we established, for j =3,4, 
$$\liminf_{x \to \infty}\liminf_{n \to \infty}\frac{I_j}{\bar\zeta(\lambda_n)}\geq 0.$$
Thus we have proved (3.6).
\hfill $\Box$
 
\begin{lem} Let $\gamma \geq 0$.  Suppose that $\zeta \in \mathcal{OS}\cap\mathcal{L}(\gamma) $. For $j=1,2$, let $\rho_j$  be distributions on $\mathbb R_+$  satisfying (3.1).
Suppose further that, for every $a \geq 0$,
$$\lim_{x \to \infty}\frac{e^{-\gamma a}\bar\rho_1(x-a)-\bar\rho_1(x)}{\bar\zeta(x)} = 0$$
and, for every $a \geq 0$,
\begin{equation}
\lim_{x \to \infty}\frac{e^{-\gamma a}\overline{\rho_1*\rho_2}(x-a)-
\overline{\rho_1*\rho_2}(x)}{\bar\zeta(x)} = 0
\end{equation}
and that $e^{\gamma x}\rho_1(dx) \in \mathcal{W}$. Then, we have, for every $a \geq 0$, 
\begin{equation}
\lim_{x \to \infty}\frac{e^{-\gamma a}\bar\rho_2(x-a)-\bar\rho_2(x)}
{\bar\zeta(x)} = 0.
\end{equation}
\end{lem}
Proof. Let $\Lambda_2$ be the totality of increasing sequences $\{\lambda_n\}_{n=1}^{\infty}$ with $\lim_{n \to\infty} \lambda_n = \infty$ such that, for every $x \in \mathbb R$, the following  limit  exists and is finite:
$$m_2(x;\{\lambda_n\}):=\lim_{n \to\infty}\frac{\bar\rho_2(\lambda_n-x)} {\bar\zeta(\lambda_n)}.$$
We have $\Lambda_2\subset \Lambda$. As for $\Lambda$, it follows that, under the assumption that $\zeta \in\mathcal{OS}$ and $\overline{\rho_2}(x) \leq C_2\bar\zeta(x)$, there exists an increasing subsequence $\{\lambda_n\} \in \Lambda_2 $ of $\{x_n\}$ for each sequence   
 $\{x_n\}_{n=1}^{\infty}$ with $\lim_{n \to\infty} x_n = \infty$. Let $\{\lambda_n\} \in \Lambda_2. $
  Recall from Lemma 2.4  that $m(x;\{\lambda_n\}) =e^{\gamma x}$ and $\widehat\zeta(\gamma) < \infty.$ As in the proof of Lemma 3.2, we have (3.5).
We find that, for every $a\in \mathbb R$,
\begin{equation}
\begin{split}
&l(x):=\lim_{n \to \infty}\frac{I_4}{\bar\zeta(\lambda_n)}\\ \nonumber
&=\int_{0-}^{x+}(e^{-\gamma a}m_2(a+y;\{\lambda_n\})-m_2(y;\{\lambda_n\}))\rho_1(dy).
\end{split}
\end{equation}
Define $M_2(y;\{\lambda_n\}):=e^{-\gamma y}m_2(y;\{\lambda_n\})$.
Then $M_2(y;\{\lambda_n\})\leq C_2$  on $\mathbb R$. Note that
$$l(x)=\int_{0-}^{x+}(M_2(a+y;\{\lambda_n\})-M_2(y;\{\lambda_n\}))e^{\gamma y}\rho_1(dy).$$
 We see from (3.2), (3.3) of Lemma 3.1, (3.5), and (3.7) that, for every $a\in \mathbb R$,
$$\lim_{x \to \infty}l(x)=\int_{0-}^{\infty}(M_2(a+y;\{\lambda_n\})-M_2(y;\{\lambda_n\}))e^{\gamma y}\rho_1(dy)=0.$$
Thus we obtain that, for every $a, b \in \mathbb R$,
$$\int_{0-}^{\infty}(M_2(a+b+y;\{\lambda_n\})-M_2(b+y;\{\lambda_n\}))e^{\gamma y}\rho_1(dy)=0.$$
Since $e^{\gamma y}\rho_1(dy) \in \mathcal{W}$, we find from Lemma 2.7 that, for every $a\in \mathbb R$,
$$M_2(a+b;\{\lambda_n\})=M_2(b;\{\lambda_n\})\mbox{ for a.e. } b \in \mathbb R.$$
Since the function $m_2(x;\{\lambda_n\})$ is increasing, the functions  $M_2(x+;\{\lambda_n\})$ and  $M_2(x-;\{\lambda_n\})$ exist for all $x \in \mathbb R$.
Taking $b_n=b_n(a)\downarrow 0$ and $b_n=b_n(a)\uparrow 0$, we have
$$M_2(a+;\{\lambda_n\})=M_2(0+;\{\lambda_n\})\mbox{ and } M_2(a-;\{\lambda_n\})=M_2(0-;\{\lambda_n\}).$$
As $a\uparrow 0$ in the first equality, we see that
$$M_2(0-;\{\lambda_n\})=M_2(0+;\{\lambda_n\})$$
and hence, for every $a\in \mathbb R$,
$$M_2(a;\{\lambda_n\})=M_2(0;\{\lambda_n\}).$$
Since $\{\lambda_n\} \in \Lambda_2 $ is arbitrary, we have (3.8).
\hfill$\Box$
\begin{lem} Let $\gamma \geq 0$. Suppose that $\zeta \in \mathcal{OS} $. For $j=1,2$, let $\rho_j$  be distributions on $\mathbb R_+$  satisfying (3.1). Suppose further that, for $j=1,2,$ and for every $a \geq 0$,
$$\liminf_{x \to \infty}\frac{e^{-\gamma a}\bar\rho_j(x-a)-\bar\rho_j(x)}{\bar\zeta(x)} \geq 0.$$
If we have, for every $a \geq 0$,
$$\lim_{x \to \infty}\frac{e^{-\gamma a}\overline{\rho_1*\rho_2}(x-a)-\overline{\rho_1*\rho_2}(x)}{\bar\zeta(x)} = 0,$$
then, for $j=1,2,$ and for every $a \geq 0$,
$$\lim_{x \to \infty}\frac{e^{-\gamma a}\bar\rho_j(x-a)-\bar\rho_j(x)}{\bar\zeta(x)} = 0.$$
\end{lem}
Proof.  Suppose that, for some $a >0$, 
$$\limsup_{x \to\infty}\frac{e^{-\gamma a}\overline{\rho_2}(x-a)-\overline{\rho_2}(x)}{\bar\zeta(x)} >0.$$
Then there is $\{\lambda_n\} \in \Lambda$ such that, for some $a >0$,
$$\lim_{n \to\infty}\frac{e^{-\gamma a}\overline{\rho_2}(\lambda_n-a)-\overline{\rho_2}(\lambda_n)}{\bar\zeta(\lambda_n)} =: \delta_0 >0.$$
So there is $\delta_1 >0$ such that, for some $a > 0$,
$$\liminf_{n \to\infty}\frac{e^{-\gamma(a+\delta_1)}\overline{\rho_2}(\lambda_n-a)-\overline{\rho_2}(\lambda_n)}{\bar\zeta(\lambda_n)} =: \delta_2>0.$$
Take $y_0$ such that $x >y_0 >\delta_1$ and $\rho_1((y_0-\delta_1,y_0]) >0$.
Let $\lambda_n':=\lambda_n+y_0$ and $a':=a+\delta_1$. Then we have
\begin{equation}
\begin{split}
&\int_{y_0-\delta_1}^{y_0}(e^{-\gamma a'}\overline{\rho_2}(\lambda_n'-a'-y)-\overline{\rho_2}(\lambda_n'-y))\rho_1(dy)\\
&\geq \rho_1((y_0-\delta_1,y_0])(e^{-\gamma  a'}\overline{\rho_2}(\lambda_n-a)-\overline{\rho_2}(\lambda_n)).
\end{split}
\end{equation}
Let $\lambda_n' > a'+x$ and $x >0$. Define $J$ as
\begin{equation}
J:=e^{-\gamma a'}\overline{\rho_1*\rho_2}(\lambda_n'-a')-\overline{\rho_1*\rho_2}(\lambda_n'). \nonumber
\end{equation}
Then we have as in  assertion (i) of Lemma 3.1
$$J=\sum_{j=1}^{4} I_j',$$
where
$$I_1':=-\int_{\lambda_n'-a'-x}^{\lambda_n'-x}\overline{\rho_1}(\lambda_n'-y)\rho_2(dy),$$
$$I_2':=\overline{\rho_1}(x)(e^{-\gamma a'}\overline{\rho_2}(\lambda_n'-a'-x)-\overline{\rho_2}(\lambda_n'-x)),$$
$$I_3':=\int_{0-}^{(\lambda_n'-a'-x)+}(e^{-\gamma a'}\overline{\rho_1}(\lambda_n'-a'-y)-\overline{\rho_1}(\lambda_n'-y))\rho_2(dy),$$
and
$$I_4':=\int_{0-}^{x+}(e^{-\gamma a'}\overline{\rho_2}(\lambda_n'-a'-y)-\overline{\rho_2}(\lambda_n'-y))\rho_1(dy).$$
For $1\le j\le 3,$ let
$$J_j:=I_j',$$
and let
$$
I_4'=\sum_{j=4}^{6}J_j,$$
where
$$J_4:=\int_{0-}^{(y_0-\delta_1)+}(e^{-\gamma a'}\overline{\rho_2}(\lambda_n'-a'-y)-\overline{\rho_2}(\lambda_n'-y))\rho_1(dy),$$
$$J_5:=\int_{y_0}^x(e^{-\gamma a'}\overline{\rho_2}(\lambda_n'-a'-y)-\overline{\rho_2}(\lambda_n'-y))\rho_1(dy),$$
and
$$J_6:=\int_{y_0-\delta_1}^{y_0}(e^{-\gamma a'}\overline{\rho_2}(\lambda_n'-a'-y)-\overline{\rho_2}(\lambda_n'-y))\rho_1(dy).$$
Then we have 
$$J=\sum_{j=1}^{6} J_j.$$
As in the proof of Lemma 3.3, we see from the assumption  and (3.9) that
\begin{equation}
\begin{split}
0&=\lim_{n\to\infty}\frac{J}{\bar\zeta(\lambda_n')}\\ \nonumber
&\geq \sum_{j=1}^{6}\liminf_{x\to\infty}\liminf_{n\to\infty}\frac{J_j}{\bar\zeta(\lambda_n')}\\
&\geq\liminf_{n\to\infty}\frac{J_6}{\bar\zeta(\lambda_n')}\\
&\geq\liminf_{n\to\infty}\rho_1((y_0-\delta_1,y_0])\frac{(e^{-\gamma a'}\overline{\rho_2}(\lambda_n-a)-\overline{\rho_2}(\lambda_n))}{{\bar\zeta(\lambda_n')}}\\
&=\rho_1((y_0-\delta_1,y_0])\frac{\delta_2}{m(-y_0;\{\lambda_n\})} >0.
\end{split}
\end{equation}
This is a contradiction. Thus we have,
 for every $a \geq 0$,
$$\lim_{x \to \infty}\frac{e^{-\gamma a}\bar\rho_2(x-a)-\bar\rho_2(x)}{\bar\zeta(x)} = 0.$$
By the analogous argument, we have
 for every $a \geq 0$,
$$\lim_{x \to \infty}\frac{e^{-\gamma a}\bar\rho_1(x-a)-\bar\rho_1(x)}{\bar\zeta(x)} = 0.$$
Thus we have proved the lemma. \hfill $\Box$
\begin{lem} Let $\gamma \geq 0$. Let $\rho$ be a distribution on $\mathbb R_+$. Suppose that  $\rho \in \mathcal{OS}$ and, for every $a \geq 0$, 
\begin{equation}
\liminf_{x \to \infty}e^{-\gamma a}\bar \rho(x-a)/\bar \rho(x) \geq 1.
\end{equation}
Then, for some positive integer $n\geq 2$,
$\rho^{n*} \in \mathcal{L}(\gamma)$ implies that $\rho \in \mathcal{L}(\gamma)$.
\end{lem}
Proof. Let $\zeta:=\rho$. Then we see from Lemma 3.3 that, for every $k \in\mathbb N$ and every  $a \geq 0$, 
$$\liminf_{x \to \infty}\frac{e^{-\gamma a}\overline{\rho^{k*}}(x-a)-\overline{\rho^{k*}}(x)}{\bar\rho(x)} \geq 0.$$
Thus we find that $\rho_1:=\rho$ and $\rho_2:=\rho^{(n-1)*}$ satisfy the assumptions of Lemma 3.5. Hence we have by Lemma 3.5,  for every $a \geq 0$,
$$\lim_{x \to \infty}\frac{e^{-\gamma a}\bar\rho(x-a)-\bar\rho(x)}{\bar\rho(x)} = 0.$$
That is, $\rho \in \mathcal{L}(\gamma)$. \hfill $\Box$
\begin{rem}
For $\gamma = 0$, the assumption (3.10) necessarily holds, but for  $\gamma > 0$, without the assumption (3.10) the lemma does not hold. For $\gamma > 0$, Watanabe \cite{w2} made a distribution $\eta \in \mathcal{OS}$ such that $\eta^{n*} \in \mathcal{L}(\gamma)$ for every $n \geq 2$ but $\eta \notin \mathcal{L}(\gamma)$.

\end{rem}
\section{Proof of results}
In this section, we prove the results stated in Sect. 1.

\begin{lem}  Let $\gamma \geq 0$ and $\mu \in \mathcal{OS}\cap\mathcal{ID}$.  If, for every $a \geq 0$, (1.1) holds,
 then, 
  for all $n\in \mathbb N$ and  every $a \geq 0$, 
\begin{equation}
\lim_{x \to \infty}\frac{e^{-\gamma a}\overline{\nu_1^{n*}}
(x-a)-\overline{\nu_1^{n*}}(x)}
{ \overline{\nu_1^{n_0*}}(x)}=0,
 \end{equation}
and we have  $T(\mu,\gamma)=(0,\infty)$.
 \end{lem}
Proof. By induction, we see from Lemma 3.2 that if (1.1) holds for every $a \geq 0$,
 then, 
  for all $n\in \mathbb N$ and  every $a \geq 0$, we have (4.1).
We have with $c:=\bar\nu(1)$, for $t>0$,
$$\mu_1^{t*}:=e^{-ct}\sum_{k=0}^{\infty}\frac{(ct)^k}{k!}\nu_1^{k*}.$$
 Suppose that, 
  for all $n\in \mathbb N$ and  every $a \geq 0$, (4.1) holds. Let $\epsilon > 0$ be arbitrary.  By Lemma 2.3,  we can choose sufficiently large $N \in \mathbb N$ such that, for  $\epsilon > 0$,
$$e^{-ct}\sum_{k=N+1}^{\infty}\frac{(ct)^k}{k!}\frac{|e^{-\gamma a}\overline{\nu_1^{k*}}
(x-a)-\overline{\nu_1^{k*}}(x)|}
{ \overline{\nu_1^{n_0*}}(x)} < \epsilon.$$
We find from (4.1) that, for every $a \geq 0$,
$$\lim_{x \to \infty}e^{-ct}\sum_{k=1}^{N}\frac{(ct)^k}{k!}\frac{e^{-\gamma a}\overline{\nu_1^{k*}}
(x-a)-\overline{\nu_1^{k*}}(x)}
{ \overline{\nu_1^{n_0*}}(x)} =0.$$
Thus we see   that, for every $a \geq 0$ and for every $t>0$,
\begin{equation}
 \lim_{x \to \infty}\frac{e^{-\gamma a}\overline{\mu_1^{t*}}(x-a)-\overline{\mu_1^{t*}}(x)}{ \overline{\nu_1^{n_0*}}(x)} =0. \nonumber
 \end{equation}
Since $\overline{\mu_1^{t*}}(x)\asymp \overline{\nu_1^{n_0*}}(x)$ for every $t>0$, we have $T(\mu,\gamma)=(0,\infty)$. \hfill $\Box$
 \begin{lem} Let $\gamma \geq 0$ and $\mu \in \mathcal{OS}\cap\mathcal{ID}$. If $0$ is a limit point of $T(\mu,\gamma)$, then, for every $a \geq 0$, (1.1) holds. 
 \end{lem}
Proof. Suppose that $0$ is a limit point of $T(\mu,\gamma)$.
Then,   there exists a strictly decreasing sequence $\{t_n\}_{n=1}^{\infty}$ in $T(\mu,\gamma)$ converging to $0$ as $n \to \infty$. We have with $c:=\bar\nu(1)$
$$\mu_1^{t_n*}:=e^{-ct_n}\sum_{k=0}^{\infty}\frac{(ct_n)^k}{k!}\nu_1^{k*}.$$
Since $\{t_n\}_{n=1}^{\infty}$ in $T(\mu,\gamma)$ and $\overline{\mu_1^{t_n*}}(x) \asymp \overline{\nu_1^{n_0*}}(x)$ from Lemma 2.2, we see that, for every $a \geq 0$,
\begin{equation}
 \begin{split}
 &\lim_{x \to \infty}\frac{e^{-\gamma a}\overline{\mu_1^{t_n*}}(x-a)-\overline{\mu_1^{t_n*}}(x)}{ \overline{\nu_1^{n_0*}}(x)}\\
 & =\lim_{x \to \infty}\frac{e^{-\gamma a}\overline{\mu_1^{t_n*}}(x-a)-\overline{\mu_1^{t_n*}}(x)}{ \overline{\mu_1^{t_n*}}(x)}\frac{ \overline{\mu_1^{t_n*}}(x)}{\overline{\nu_1^{n_0*}}(x)}=0. \nonumber
\end{split} 
\end{equation}
Thus we obtain from Lemma 2.3  that, for every $a \geq 0$,
\begin{equation}
\begin{split}
&\limsup_{x \to \infty}|\frac{e^{-\gamma a}\overline{\nu_1}(x-a)-\overline{\nu_1}(x)}
{ \overline{\nu_1^{n_0*}}(x)}| \\
&=\limsup_{n \to \infty}\limsup_{x \to \infty}|\frac{e^{ct_n}}{ct_n}\frac{e^{-\gamma a}\overline{\mu_1^{t_n*}}(x-a)-\overline{\mu_1^{t_n*}}(x)}{ \overline{\nu_1^{n_0*}}(x)}
-\frac{e^{-\gamma a}\overline{\nu_1}(x-a)-\overline{\nu_1}(x)}
{ \overline{\nu_1^{n_0*}}(x)}|\\
&\leq\limsup_{n \to \infty}\limsup_{x \to \infty}\sum_{k=2}^{\infty}\frac{(ct_n)^{(k-1)}}{k!}\frac{e^{-\gamma a}\overline{\nu_1^{k*}}(x-a)+\overline{\nu_1^{k*}}(x)}{ \overline{\nu_1^{n_0*}}(x)} =0. \nonumber
\end{split}
 \end{equation}
Thus we have (1.1) for every $a \geq 0$.
\hfill $\Box$

\begin{lem} Let $\gamma \geq 0$ and $\mu \in \mathcal{OS}\cap\mathcal{ID}$. If $t_0, t_1 \in T(\mu,\gamma)$ with $t_1 >t_0$, then $t_1-t_0 \in T(\mu,\gamma)$. If $T(\mu,\gamma)$ has a limit point, then $T(\mu,\gamma) =(0,\infty)$. If $T(\mu,\gamma)$ has the minimum $a_0 >0$, then
$T(\mu,\gamma) =a_0\mathbb N$.
\end{lem}
Proof. Suppose that $t_0, t_1 \in T(\mu,\gamma)$ with $t_1 >t_0$. Let $\zeta:=\rho_1:=\mu^{t_0*}$ and $\rho_2:=\mu^{(t_1-t_0)*}$.   The distribution $e^{\gamma x}\rho_1(dx)/\widehat\rho_1(\gamma)$ is an exponentially tilted infinitely divisible distribution and hence itself is infinitely divisible, thus having a non-vanishing characteristic function. That is, $e^{\gamma x}\rho_1(dx) \in \mathcal{W}$. See (iii) of Theorem 25.17  of Sato \cite{s}. Thus we see from Lemma 3.4 that $\mu^{(t_1-t_0)*}\in \mathcal{L}(\gamma)$. Thus, if $T(\mu,\gamma)$ has a limit point, then 0 is a limit point of $T(\mu,\gamma)$, and hence, by Lemmas 4.1 and 4.2, $T(\mu,\gamma) =(0,\infty)$. If $T(\mu,\gamma)$ has the minimum $a_0 >0$, then clearly  $a_0\mathbb N\subset T(\mu,\gamma)$ and $ T(\mu,\gamma)\setminus a_0\mathbb N =\emptyset $.  \hfill$\Box$

Proof of Theorem 1.1. Assertion (i) is clear from Lemmas 4.1, 4.2, and 4.3. The first part of assertion (ii) is due to Lemmas 4.1   and 4.2. Suppose that $\mu \in \mathcal{A}(\gamma)$. If $n < n_0,$ then $\nu_1^{n*} \not\in \mathcal{OS}$ simply because of the definition of $n_0$. If $n \ge n_0$ and $x$ is large, then $\overline{\nu_1^{n*}}(x) \ge \overline{\nu_1^{n_0*}}(x)$, and hence (4.1) implies that $\nu_1^{n*} \in \mathcal{L}(\gamma)$.
\hfill$\Box$

Proof of Corollary 1.1. Suppose that $\mathcal{C}(\gamma)$ is not empty.  Then there is the minimum  $a_0 >0$ in $T(\mu,\gamma)$ for $\mu \in \mathcal{C}(\gamma)$. Since $a_0>0$ is a period of $T(\mu,\gamma)$, for  $n=2$, $\mu^{a_0*}=(\mu^{(a_0/n)*})^{n*} \in \mathcal{L}(\gamma)$ but $(\mu^{(a_0/n)*})^{(n+1)*} \notin \mathcal{L}(\gamma)$ and $\mu^{(a_0/n)*} \notin \mathcal{L}(\gamma)$. Thus Hypotheses I and II are not true. Suppose that $\mathcal{C}(\gamma)$ is empty. Then, obviously, Hypotheses I and II are  true. Thus (1), (2), and (3) are equivalent.
We prove the equivalence of (3) and (4). Suppose that $\mathcal{C}(\gamma)$ is empty. Then for every $\mu \in \mathcal{OS}\cap\mathcal{ID}$ it holds that, for every $2t \in T(\mu,\gamma)$, $\mu_1^{t*}\in \mathcal{L}(\gamma)$ and hence, for  all $a \geq 0$, (1.2) holds.
Conversely, suppose that $\mathcal{C}(\gamma)$ is not empty and, for  $a_0=2t \in T(\mu,\gamma)$ with $\mu \in\mathcal{C}(\gamma)$ and for  all $a \geq 0$, (1.2) holds.
Letting  $\rho_1:=\rho_2:=\mu_1^{t*}$, $\zeta:=\mu_1^{2t*}$, define $\Lambda_2$ as in Lemma 3.4 and let $\{\lambda_n\} \in \Lambda_2\subset \Lambda$. We have (3.3) by Lemma 3.1
 for $j=1,2$.
We have
$I_3+I_4=2I_4+I_5$, where 
$$I_5:=\int_x^{\lambda_n-a-x}(e^{-\gamma a}\overline{\rho_1}(\lambda_n-a-y)-\overline{\rho_1}(\lambda_n-y))\rho_2(dy),$$
We have by the assumption (1.2) for every $a \geq 0$
\begin{equation}
\limsup_{x\to \infty}\limsup_{n \to\infty}\frac{|I_5|}{\bar\zeta(\lambda_n)}=0.\nonumber
\end{equation}
Define $M_2(y;\{\lambda_n\}):=e^{-\gamma y}m_2(y;\{\lambda_n\})$. 
Thus we find from (3.2), (3.3), and $2t \in T(\mu,\gamma)$ that, for every $a \geq 0$,
\begin{equation}
\begin{split}
&\lim_{x\to \infty}\lim_{n \to \infty}\frac{I_4}{\bar\zeta(\lambda_n)}\\ \nonumber
&=\int_{0-}^{\infty}(e^{-\gamma a}m_2(a+y;\{\lambda_n\})-m_2(y;\{\lambda_n\}))\rho_1(dy)\\
&=\int_{0-}^{\infty}(M_2(a+y;\{\lambda_n\})-M_2(y;\{\lambda_n\}))e^{\gamma y}\rho_1(dy)=0.
\end{split}
\end{equation}
 The distribution $e^{\gamma x}\rho_1(dx)/\widehat\rho_1(\gamma)$ is an exponentially tilted infinitely divisible distribution and hence itself is infinitely divisible, thus having a non-vanishing characteristic function. That is, 
$$e^{\gamma y}\rho_1(dy)=e^{\gamma y}\mu_1^{t*}(dy)\in \mathcal{W}.$$
As in the proof of Lemma 3.4, we have $\rho_2=\mu_1^{t*}\in \mathcal{L}(\gamma)$. This is a contradiction. Thus (3) and (4) are equivalent.
\hfill $\Box$

Proof of Remark 1.1. Let $\gamma =0$. Then we see from Lemma 3.6 that Hypothesis II is true. Thus $\mathcal{C}(0)$ is empty and hence Remark 1.1 follows from Theorem 1.1.  \hfill$\Box$

Proof of Proposition 1.1.  Let $\gamma > 0$ and $\mu \in \mathcal{OS}\cap\mathcal{ID}$.  Suppose that (1.3) holds for every $a \geq 0$.
Let $\zeta :=\nu_1^{n_0*}$. Then, by induction, we see from (1.3) and Lemma 3.3 that, for every $n \in \mathbb N$ and every $a \geq 0$,
$$\liminf_{x \to \infty}\frac{e^{-\gamma a}\overline{\nu_1^{n*}}(x-a)-\overline{\nu_1^{n*}}(x)}{\overline{\nu_1^{n_0*}}(x)} \geq 0.$$
Let $\epsilon >0$ be arbitrary. Thus, letting $N \in\mathbb N$ sufficiently large, we have, for every $t >0$ and for every $a \geq 0$,
\begin{equation}
\begin{split}
&\liminf_{x \to \infty}\frac{e^{-\gamma a}\overline{\mu_1^{t*}}(x-a)-\overline{\mu_1^{t*}}(x)}
{ \overline{\nu_1^{n_0*}}(x)}\\
&=\liminf_{x \to \infty}e^{-ct}\sum_{k=1}^N\frac{(ct)^k}{k!}\frac{e^{-\gamma a}\overline{\nu_1^{k*}}
(x-a)-\overline{\nu_1^{k*}}(x)}
{ \overline{\nu_1^{n_0*}}(x)} \\
&- \limsup_{x \to \infty}e^{-ct}\sum_{k=N+1}^{\infty}\frac{(ct)^k}{k!}\frac{e^{-\gamma a}\overline{\nu_1^{k*}}
(x-a)+\overline{\nu_1^{k*}}(x)}
{ \overline{\nu_1^{n_0*}}(x)}\geq -\epsilon.\nonumber
\end{split}
 \end{equation}
Since $\epsilon >0$ is arbitrary and $\overline{\nu_1^{n_0*}}(x)\asymp \overline{\mu_1^{(t/n)*}}(x)$ for every $ n \in \mathbb N$, we obtain that $\rho:=\mu_1^{(t/n)*}$ satisfies $\rho \in \mathcal{OS}$ and (3.10) holds. Hence we find from Lemma 3.6 that if $t \in T(\mu,\gamma)$, then $t/n \in T(\mu,\gamma)$ for every $ n \in \mathbb N$. Thus, by Lemmas 4.1 and 4.2, either $T(\mu,\gamma)= (0,\infty)$ or $\emptyset$.   \hfill  $\Box$

Proof of Proposition 1.2.  Suppose that $\nu_1^{2*} \in \mathcal{L}(\gamma)$ and the real part of $\widehat\nu_1(\gamma+iz)$ is not $0$ for every $z \in\mathbb R$. If $t \in T(\mu,\gamma)$, then
\begin{equation}
\mu_1^{t*}=e^{-ct}\sum_{k=0}^{\infty}\frac{(ct)^k}{k!}\nu_1^{k*}\in \mathcal{L}(\gamma)\cap\mathcal{OS}.
\end{equation}
Define distributions $\eta_1$ and $\eta_2$ on $\mathbb R_+$ as
$$\eta_1:=(\cosh(ct))^{-1}\sum_{k=0}^{\infty}\frac{(ct)^{2k}}{(2k)!}\nu_1^{(2k)*}$$
and
$$\eta_2:=(\sinh(ct))^{-1}\sum_{k=0}^{\infty}\frac{(ct)^{2k+1}}{(2k+1)!}\nu_1^{(2k+1)*}.$$
We see from  Proposition 3.1 of Shimura and Watanabe \cite{sw1} that $\eta_j\in \mathcal{OS}$ and $\overline{\eta_j}(x)\asymp\overline{\nu_1^{n_0*}}(x)$ for $j=1,2$. Let $\epsilon >0$ be arbitrary. We obtain from Lemma 2.3 that there is a positive integer $N=N(a,\epsilon,t)$
such that
$$\limsup_{x \to\infty}(\cosh(ct))^{-1}\sum_{k=N+1}^{\infty}\frac{(ct)^{2k}}{(2k)!}\frac{e^{-\gamma a}\overline{\nu_1^{(2k)*}}(x-a)+\overline{\nu_1^{(2k)*}}(x)}{\overline{\nu_1^{n_0*}}(x)}<\epsilon.$$
Since $\nu_1^{(2k)*}\in \mathcal{L}(\gamma)$ for every $k \ge 0$, we have, for every $a \ge 0$ and every $t >0$,
$$\limsup_{x \to\infty}(\cosh(ct))^{-1}\sum_{k=0}^{N}\frac{(ct)^{2k}}{(2k)!}\frac{|e^{-\gamma a}\overline{\nu_1^{(2k)*}}(x-a)-\overline{\nu_1^{(2k)*}}(x)|}{\overline{\nu_1^{n_0*}}(x)}=0.$$
Thus with some $C=C(t) >0$ we have, for every $a \ge 0$ and every $t >0$,
\begin{equation}
\begin{split}
&\limsup_{x \to\infty}\frac{|e^{-\gamma a}\overline{\eta_1}(x-a)-\overline{\eta_1}(x)|}{\overline{\eta_1}(x)} \nonumber \\
&\le \limsup_{x \to\infty}(\cosh(ct))^{-1}\sum_{k=0}^{N}\frac{(ct)^{2k}}{(2k)!}\frac{|e^{-\gamma a}\overline{\nu_1^{(2k)*}}(x-a)-\overline{\nu_1^{(2k)*}}(x)|}{C\overline{\nu_1^{n_0*}}(x)}\\
&+\limsup_{x \to\infty}(\cosh(ct))^{-1}\sum_{k=N+1}^{\infty}\frac{(ct)^{2k}}{(2k)!}\frac{e^{-\gamma a}\overline{\nu_1^{(2k)*}}(x-a)+\overline{\nu_1^{(2k)*}}(x)}{C\overline{\nu_1^{n_0*}}(x)}\\
& \le \epsilon/C.
\end{split}
\end{equation}
Since $\epsilon>0$ is arbitrary, we have
\begin{equation}
\eta_1\in \mathcal{L}(\gamma)\cap\mathcal{OS}.
\end{equation}
Since
$$\sinh(ct)\eta_2=e^{ct}\mu_1^{t*} -\cosh(ct)\eta_1,$$
we have by (4.2) and (4.3)
$$\eta_2\in \mathcal{L}(\gamma)\cap\mathcal{OS}.$$
Let $\zeta:=\rho_1:=\eta_2$ and $\rho_2:=\nu_1.$ Then, by argument similar to the proof of  (4.3),
$$\rho_1*\rho_2=(\sinh(ct))^{-1}\sum_{k=0}^{\infty}\frac{(ct)^{2k+1}}{(2k+1)!}\nu_1^{(2k+2)*}\in \mathcal{L}(\gamma)\cap\mathcal{OS}.$$
Since the real part of $\widehat\nu_1(\gamma+iz)$ is not $0$ for every $z \in\mathbb R$,
$$2\sinh(ct)\widehat\rho_1(\gamma+iz)=\exp(ct\widehat\nu_1(\gamma+iz))-\exp(-ct\widehat\nu_1(\gamma+iz))\ne 0$$ 
for every $z \in\mathbb R$, that is, $e^{\gamma x}\rho_1(dx) \in \mathcal{W}$. Thus we see from Lemma 3.4 that 
$$\lim_{x \to \infty}\frac{e^{-\gamma a}\bar\nu_1(x-a)-\bar\nu_1(x)}
{\bar\zeta(x)} = 0.$$
Since $\bar\zeta(x)\asymp\overline{\nu_1^{n_0*}}(x)$, we see from Theorem 1.1 that $T(\mu,\gamma)=(0,\infty)$. Thus we have proved the proposition.  \hfill$\Box$

Proof of Proposition 1.3.  Let $\gamma > 0$ and $\mu \in \mathcal{OS}\cap\mathcal{ID}$.  Suppose that $\nu_1^{n_1*} \in \mathcal{S}_{\sharp}$.  Since $\overline{\mu^{t*}}(x)\asymp \overline{\nu_1^{n_1*} }(x)$, we have
$\mu^{t*}\in \mathcal{S}_{\sharp}$ for every $t >0$. Thus we see from Lemmas 2.5 and 2.6 that if  $ T(\mu,\gamma)\ne \emptyset$, then $\nu_1 \in \mathcal{S}(\gamma)$ and hence $T(\mu,\gamma)=(0,\infty)$. That is, either $T(\mu,\gamma)= (0,\infty)$ or $\emptyset$. Moreover, $T(\mu,\gamma)=(0,\infty)$ if and only if $\nu_1 \in \mathcal{S}(\gamma)$.  \hfill$\Box$

\medskip

{\bf Acknowledgements}

The author is grateful to the referee for his careful reading the manuscript and helpful advice.


\medskip
\begin{thebibliography}{99}
\bibitem{cwx} Cui, Z.,  Wang, Y., Xu, H. : Some positive conclusions related to the Embrechts-Goldie conjecture. arXiv:1609.00912 (2016).
\bibitem{eg1}  Embrechts, P., Goldie, C. M. : On closure and factorization properties of subexponential and related distributions. J. Aust. Math. Soc. A 29 (1980) 243–256. 
\bibitem{eg2} Embrechts, P., Goldie, C. M. : On convolution tails. Stochastic Process. Appl. 13 (1982) 263–278. 
 \bibitem{egv} Embrechts, P., Goldie, C.M., Veraverbeke, N.: Subexponentiality and infinite divisibility.  Z. Wahrscheinlichkeitstheorie Verw. Gebiet. 49 (1979) 335-347.
\bibitem{k} Kl\"uppelberg, C. :  Asymptotic ordering of distribution functions on convolution semigroup. Semigroup Forum. 40 (1990) 77–92. 
 \bibitem{kv} Kl\"uppelberg, C., Villasenor, J. A. : The full solution of the convolution closure problem for convolution-equivalent distributions. J. Math. Anal. Appl. 160 (1991) 79-92. 
 \bibitem{l} Leslie, J. R. : On the non-closure under convolution of the subexponential family. J. Appl. Probab. 26 (1989) 58-66.
 \bibitem{p1} Pakes, A.G.: Convolution equivalence and infinite divisibility. J. Appl. Probab. 41 (2004) 407-424.
\bibitem{p2} Pakes, A. G. : Convolution equivalence and infinite divisibility: Corrections and corollaries. J. Appl. Probab. 44 (2007) 295–305. 
\bibitem{s} Sato, K.: L\'evy processes and infinitely divisible distributions. Cambridge Studies in Advanced Mathematics, 68 Cambridge Univ. Press. 2013. 
\bibitem{sw1} Shimura, T., Watanabe, T.: Infinite divisibility and generalized subexponentiality. Bernoulli 11 (2005) 445-469. 
\bibitem{sw2} Shimura, T., Watanabe, T. :  On the convolution roots in the convolutionequivalent class. The Institute of Statistical Mathematics Cooperative Research Report 175 (2005) pp. 1–15. 
\bibitem{w1} Watanabe, T.: Convolution equivalence and distributions of random sums. Probab. Theory Related Fields 142 (2008) 367-397. 
\bibitem{w2} Watanabe, T.: The Wiener condition and the conjectures of Embrechts and Goldie. Ann. Probab. 47 (2019) 1221–1239. 
\bibitem{wy1} Watanabe, T., Yamamuro, K.:  Ratio of the tail of an infinitely divisible distribution on the line to that of its 
L\'evy measure. Electron. J. Probab. 15 (2010) 44-74.
\bibitem{wy2} Watanabe, T., Yamamuro, K.: Two non-closure properties on the class of subexponential densities. J. Theoret. Probab. 30 (2017) 1059-1075.
\bibitem{w}  Wiener, N.:  Tauberian theorems. Ann. of Math. (2)33 (1932) 1–100.
\bibitem{xfw} Xu, H., Foss, S., Wang, Y. : Convolution and convolution-root properties of longtailed distributions. Extremes 18 (2015) 605–628. 
\bibitem{xwcy} Xu, H., Wang, Y.,  Cheng, D.,  Yu, C.: On the closure under convolution roots related to an infinitely divisible distribution in the distribution class $\mathcal{L}(\gamma)$. arXiv:1512.01792 (2015).
\end{thebibliography}
\end{document}